\documentclass{amsart} 

\usepackage{amsmath,amssymb}
\usepackage{amsthm}
\usepackage[all]{xy}

\newcommand{\e}[2]{{\tilde e}(#1,#2)}
\def\Xt{{\widetilde X}}
\def\Xb{{\bar X}}
\def\Cn{C_n}
\def\Cnb{{\bar C_n}}
\def\Cnt{{\widetilde C_n}}
\def\Cnh{{\hat C_n}}
\def\ra{\rightarrow}
\def\K{\mathcal{K}}
\def\L{\mathcal{L}}
\def\pcn{p^{-1}(C_n)}
\def\pdn{p^{-1}(D_n)}
\def\qcn{q^{-1}({\bar C_n})}
\def\cc{\underset c -}
\def\ni{\noindent}

\newtheorem{thm}{Theorem}
\newtheorem{cor}[thm]{Corollary}
\newtheorem{lem}[thm]{Lemma}
\newtheorem{rem}[thm]{Remark}
\newtheorem{claim}[thm]{Claim}

\theoremstyle{definition}
\newtheorem{defn}[thm]{Definition}

\begin{document}

\title{Filtered Ends of Infinite Covers and Groups} 
\thanks{To appear in the Journal of Pure and Applied Algebra}

\author{Tom Klein}

\begin{abstract}
  Let $f:A\ra B$ be a covering map.  We say $A$ has $e$ filtered ends with
  respect to $f$ (or $B$) if for some filtration $\{K_n\}$ of $B$ by compact
  subsets, $A - f^{-1}(K_n)$ ``eventually'' has $e$ components.  The main
  theorem states that if $Y$ is a (suitable) free $H$-space, if $K < H$ has
  infinite index, and if $Y$ has a positive finite number of filtered ends with
  respect to $H\backslash Y$, then $Y$ has one filtered end with respect to
  $K\backslash Y$.  This implies that if $G$ is a finitely generated group and
  $K < H < G$ are subgroups each having infinite index in the next, then $0 <
  \e{G}{H} < \infty$ implies $\e{G}{K} = 1$, where $\e{\cdot}{\cdot}$ is the
  number of filtered ends of a pair of groups in the sense of Kropholler and
  Roller.
\end{abstract}

\maketitle

\section{Introduction}
\label{intro}

The number of relative ends of a group $G$ and subgroup $H$, denoted
$\e{G}{H}$, was originally introduced algebraically by Kropholler and Roller
\cite{KR}.  Geoghegan \cite{Ge} gives a topological description of $\e{G}{H}$
(when $G$ is finitely generated) as the number of filtered ends of a certain
filtration (derived from $H$) of the universal cover of a presentation
complex for $G$.  To be more precise (see also section~\ref{defns}), let
$f:A\rightarrow B$ be a covering map of CW complexes of locally finite type.
We say $A$ has $e=\e{A}{f}$ filtered ends with respect to $f$ (or $B$) if for
some filtration $\{C_n\}$ of $B$ by finite subcomplexes, the cardinality of
$\underset n \varprojlim\left\{\pi_0\left(A\cc f^{-1}(C_n)\right)\right\}$ is
$e$ (this is independent of the filtration).  $\e{G}{H}$, also called the
number of filtered ends of the pair, is then $\e{\Xt^1}{q}$, where $X$ is a
path connected CW complex with $\pi_1(X)=G$, $\Xt$ is the universal cover of
$X$, and $q:\Xt^1\rightarrow \Xb(H)^1$ is the restriction to $1$-skeleta of the
natural map from $\Xt$ to the covering space of $X$ corresponding to $H$.  The
following theorem on filtered ends of covering spaces (proved in section
\ref{mtproof}) and its corollary on filtered ends of pairs of groups (proved in
section \ref{corollary}) are the main results of the article.

\begin{thm}
\label{mt}
Let $H$ be a group and let $Y$ be a path connected free $H$-CW complex of
locally finite type.  Let $K<H$ with $[H:K]=\infty$, and let $p:Y\ra
H\backslash Y$, $q:Y\ra K\backslash Y$ be the quotient covering maps.  Then $0
< \e{Y}{p} < \infty$ implies $\e{Y}{q}=1$.
\end{thm}

Applied to $1$-skeleta with $Y$ equal to the universal cover of a presentation
complex (with finite $1$-skeleton) of a finitely generated group $G$, we get
the following corollary for filtered ends of pairs of groups:

\begin{cor}
\label{cor1}
Let $G$ be a finitely generated group with $K \leq H < G$ and assume $\e{G}{H}
= n$ is finite and non-zero.  Then
\begin{displaymath} \e{G}{K} = \left\{
\begin{array}{rl}
       n & \mbox{ if } [H:K]<\infty \\
       1 & \mbox{ if } [H:K]=\infty.
\end{array}
\right.
\end{displaymath}
\end{cor}

Note also the contrapositive of Corollary \ref{cor1}: if $\e{G}{K} > 1$ and
$[G:H] = [H:K] =\infty$, then $\e{G}{H}=\infty$.  This contrapositive combined
with other results can be used to show that Thompson's groups $T$ and $V$ are
not K\"ahler (cf. \cite{NR}, Lemma 0.8 and Theorem 0.3): $T$ and $V$ are known
to have subgroups with at least two filtered ends and then the contrapositive
provides subgroups with at least three filtered ends so that the theorems of
\cite{NR} can be applied.  It seems likely that similar arguments could be
applied to other groups with subgroups having two filtered ends.

\section*{Acknowledgements} Theorem \ref{mt} is in answer to a question
of Mohan Ramachandran.  I would like to thank Ross Geoghegan for suggestions
which improved the exposition.

\section{Definitions}
\label{defns}
We now briefly review the necessary definitions (sections~4.7 and 4.9 of
\cite{Ge}).  We are heading for ``the number of ends of a CW complex with
respect to a filtration'' (Definition \ref{fe}).  We work with closed cells
and use the term {\it graph} for a $1$-dimensional CW complex. Let $X$ be an
arbitrary CW complex.

\ni $X$ is {\em locally finite} if a given cell intersects non-trivially with
only finitely many other cells, and $X$ has {\em locally finite type} if the
$n$-skeleton $X^n$ of $X$ is locally finite for all $n$. \\
\ni The {\em carrier} of a cell $e$ of X, denoted $C(e)$, is the
intersection of all subcomplexes of $X$ which contain $e$. \\
\ni The {\em CW neighborhood} of a subcomplex $A$ of $X$, denoted $N(A)$, is
the union of all cell carriers that meet $A$ non-trivially. \\
\ni The {\em CW complement} of a subcomplex $A$ of $X$, denoted $X\cc A$, is
the largest subcomplex of $X$ with $0$-skeleton $X^0 - A^0$. \\
\ni All filtrations of a CW complex $X$ will be by subcomplexes of $X$ indexed
by the natural numbers.  A {\em finite filtration} is a filtration by finite
subcomplexes\footnote{This differs from \cite{Ge}, section~3.6, which requires
  a finite filtration to be by {\it full} finite subcomplexes.}.  If $X$ has
locally finite type and $\K=\{K_n\}$ is a filtration of $X$, we say $(X,\K)$ is
{\em well filtered} if for each $n$ and $i$ there exists $j$ such that
$N_{X^n}(K_i^n) \subseteq K_j$.

Now let $(Y,\{L_i\})$ be a well filtered path connected CW complex of locally
finite type.  The {\em number of filtered ends} of $(Y,\{L_i\})$ is the
cardinality of the set $\underset i \varprojlim\left\{\pi_0\left(Y\cc
  L_i\right)\right\}$.  (We are skipping the intermediate notion of {\em
  filtered end}, which is what the number of filtered ends counts -- briefly, a
filtered end of $Y$ is an equivalence class of {\em filtered rays} $\omega :
[0,\infty) \ra Y$, where two filtered rays are equivalent if their restrictions
to $\mathbb{N}$ are {\em filtered homotopic}.  Then each $[\omega]$ picks out
an element of the inverse limit, and this correspondence is a bijection.  See
\cite{Ge} section~4.7 for more details.)
 
\begin{defn}
\label{fe}  
Let $Y$ and $X$ be path connected CW complexes of locally finite type and let
$p:Y \rightarrow X$ be a covering map.  We say {\it $Y$ has $e$ filtered ends
  over $p$} ($0 \leq e \leq \infty$), and write $\e{Y}{p}=e$, if the pair
$(Y,\{\pcn\})$ has $e$ filtered ends for some finite filtration $\{C_n\}$ of
$X$.  [By Lemma \ref{lft}.i) below there exists a finite filtration and
$(Y,\{\pcn\})$ is well filtered.  Also it is known that in this case $Y$ has
$e$ filtered ends with respect to every finite filtration of $X$.]\\
\end{defn}

\begin{rem}
\label{graph}
As noted in section~4.7 of \cite{Ge}, if $(Y,\L)$ is a well filtered path
connected CW complex of locally finite type then the number of filtered ends
of $(Y^1, \{L_i\cap Y^1\})$ equals the number of filtered ends of
$(Y,\{L_i\})$ (by the ``CW-filtered cellular approximation theorem'', \cite{Ge}
3.12.3). 
\end{rem}
We use this in the proof of Theorem \ref{mt} to reduce to the case where
$Y$ is a graph. 

Let $A$ be a subgraph of a path connected locally finite graph $X$.  In this
case $N(A)$ is simply $A \cup\,\cup\{$edges of $X$ with at least one endpoint
in $A\}$.  Then $X = N(A) \;\cup\; (X\cc A) = A \;\cup\; B \;\cup\;
(X\cc A) $, where $B = \cup\{$edges not contained in $A$ or $X\cc A\}$.  We
refer to the edges of $B$ as {\em bridge edges} (of $A$).

If $(Y,\L)$ is a filtered CW complex we say $S\subseteq Y$ is {\em
  $\L$-bounded} if $S\subseteq L_j$ for some $j$, and $(Y,\L)$ is {\em regular}
if for each $i$ the union of all $\L$-bounded path components of $Y\cc L_i$ is
$\L$-bounded (\cite{Ge}, section~4.7).

The next proposition shows that these notions all work well with respect to
covering spaces of CW complexes of locally finite type.
\begin{lem}
\label{lft}
i) If $X$ is a path connected CW complex of locally finite type then there
exists a finite filtration $\{K_n\}$, and $(X,\{K_n\})$ is well filtered (for
any finite filtration). \\
ii) If $p:Y\ra X$ is a covering map of CW complexes then $Y$ has
locally finite type if and only if $X$ does.  \\
iii) Let $p:Y\ra X$ be a covering map of CW complexes of locally finite type.
If $(X,\{K_n\})$ is well filtered then so
is $(Y,\{p^{-1}(K_n)\})$.  \\
iv) Let X be a locally finite graph with $p:Y\ra X$ a covering map, and let
$\{\Cn\}$ be a finite filtration of $X$.  Then $(Y, \{\pcn\})$ is well filtered
and regular.
\end{lem}
\begin{proof}  
  i) By \cite{Ge} 3.6.3 $X$ is countable.  Well order the cells as $e_1,
  e_2,\ldots$ and let $E_n = e_1 \cup \cdots \cup e_n$.  Since $E_n$ is compact
  it is contained in a (minimal) finite subcomplex, say $K_n$.  Then $\{K_n\}$
  is a finite filtration of $X$.  By \cite{Ge}, 3.1.12 and 3.6.9,
  $N_{X^n}(K_i^n)$ is finite, so $(X,\{K_n\})$ is well filtered. \\
  ii) $p|_{Y^n}:Y^n \ra X^n$ is a covering map, so $Y^n$ is locally finite if
  and only if $X^n$ is locally finite. \\
  iii) Let $q=p|_{Y^n}:Y^n \ra X^n$ and let $L_i =
  p^{-1}(K_i)$.\;\;\;$q(N_{Y^n}(L_i^n)) = N_{X^n}(K_i^n)$ is contained in some
  $K_j$ ($j>i$). \\ 
  iv) Well filtered follows from i) and iii).  For regular, note that $p(Y\cc
  \pcn) = X\cc \Cn$ and that if $A$ is an $\L$-bounded path component of $Y\cc
  \pcn$ ($\L=\{\pcn\}$) then $p(A)$ is contained in a finite component of $X\cc
  \Cn$ since otherwise we could construct an unbounded ray in $X\cc \Cn$
  starting in $p(A)$ which would lift to an $\L$-unbounded ray in $Y\cc \pcn$
  starting in $A$, a contradiction.  But since $N(\Cn)$ is finite and $X$ is
  locally finite, there are only finitely many finite components of $X\cc
  \Cn$.
\end{proof}

We now turn to filtered ends of pairs of groups (\cite{Ge}, section~4.9).  Let
$G$ be a finitely generated group and $H\leq G$.  Let $X$ be a path connected
CW complex with fundamental group isomorphic to $G$ and having finite
$1$-skeleton.  Let $\Xb(H)$ be the covering space of $X$ corresponding to $H$
and let $\Xt$ be the universal cover of $X$.  We work on the $1$-skeleta of
$\Xt$ and $\Xb(H)$, both locally finite graphs.  Let $p:\Xt^1 \ra \Xb(H)^1$ be
the covering projection.  The {\it number of filtered ends} of $(G,H)$ is
$\e{G}{H} := \e{\Xt^1}{p}$ (see Definition \ref{fe}).  In other words,
start with any finite filtration $\{\Cn\}$ of $\Xb(H)^1$, then count the
number of $\{\pcn\}$-unbounded components of $\Xt^1\cc p^{-1}(C_k)$ and take
the sup of that count as $k$ goes to infinity.

\section{Proof of Theorem \ref{mt}}
\label{mtproof}

\begin{rem} The case $\e{Y}{p} =1$ implies $\e{Y}{q} = 1$ of the theorem
  follows from the following ``monotonicity'' property of $\e{Y}{\cdot}$: if
  the covering map $p$ (of CW complexes of locally finite type) factors
  through the covering map $q$ as $p=r\circ q$ then $\e{Y}{q} \leq \e{Y}{p}$.
  However, we will not need this fact.
\end{rem}

{\bf PROOF of Theorem 1.}
By Remark \ref{graph} we may assume $Y$ is a graph, which we do. \\
Define $X := H\backslash Y$ and $\Xb := K\backslash Y$, and let $r:\Xb\ra X$ be
the induced covering map.  The first step of the proof is to produce compatible
filtrations of $\Xb$ and $X$ along with a rough ``fundamental domain'' in $Y$
(called $\Cnt$ below) for the $H$ and $K$ actions.  The specifics of the
filtrations are stated in the following lemma, whose proof will be postponed to
section \ref{filts}.

\begin{lem}
\label{compat} 
There exists a finite filtration $\{\Cnb\}$ of $\Xb$ and a finite filtration
$\{\Cn\}$ of $X$ with $r(\Cnb)=\Cn$ and $U_n := Y \cc \pcn$ the disjoint union
of $e$ $\{\pcn\}$-unbounded path components.  There are also finite subgraphs
$\Cnt\subseteq \pcn$ and ${\hat \Cn} \subseteq N(\pcn)$ in $Y$ with ${\hat
  \Cn}\cap \pcn = \Cnt$, $q(\Cnt) = \Cnb$, $p(\Cnt) = \Cn$, and ${\hat \Cn}
\cup U_n$ path connected.
\end{lem}
Let $\K$ be the filtration $\{\pcn\}$ of $Y$ and let $\L$ be the filtration
$\{\qcn\}$.  Define $U^1, U^2, \ldots, U^e$ to be the $e$ $\K$-unbounded path
components of $U_n = Y \cc \pcn$, so that $U_n = U^1 \sqcup \cdots \sqcup
U^e$.

$H$ preserves $\pcn$, so also bridge edges of $\pcn$ and $U_n = Y \cc \pcn$.
Then $\Cnh \cup U_n$ path connected implies \\
$\mbox{\ \ \ }(\dagger)\;\;$ $h({\hat \Cn} \cup U_n) = h{\hat \Cn} \cup U_n$ is
path connected for all $h \in H$ and all $n$.

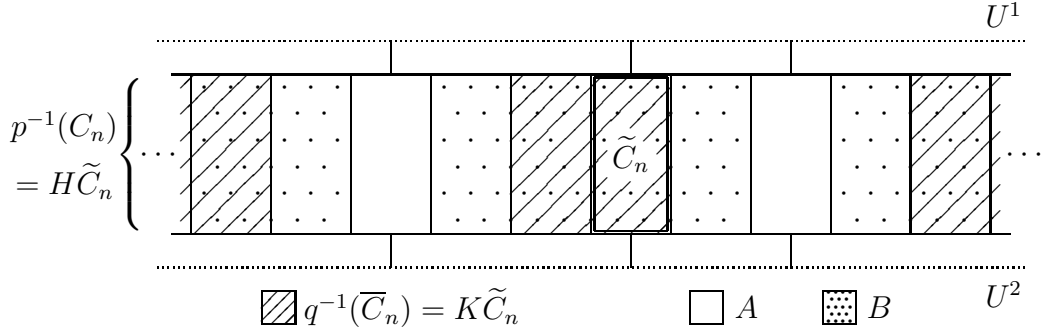
\begin{figure}[h]
{\fontsize{12}{12}\selectfont 
\def\el{2.5} 
\begin{xy}0;/r.21pc/:
(10,5); (140,5) **\dir{.};
(10,-29); (140,-29) **\dir{.};
(10,-12)*{\cdots};
(140,-12)*{\cdots};
(6,-12)*{\left\{
  \begin{xy}
  (0,12)*{ };
  (0,-12)*{ };
  \end{xy}
  \right.
};
(12,0); (138,0) **\dir{-}; (12,-24); (138,-24) **\dir{-};
(15,0); (15,-24) **\dir{-};(135,-24); (135,0) **\dir{-};
(27,0); (27,-24) **\dir{-};(39,0); (39,-24) **\dir{-};
(51,0); (51,-24) **\dir{-};(63,0); (63,-24) **\dir{-};
(75,0); (75,-24) **\dir{-};(87,0); (87,-24) **\dir{-};
(75.5,-.5); (75.5,-23.5) **\dir{-};(86.5,-.5); (86.5,-23.5) **\dir{-};
(75.5,-.5); (86.5,-.5) **\dir{-};(75.5,-23.5); (86.5,-23.5) **\dir{-};
(99,0); (99,-24) **\dir{-};(111,0); (111,-24) **\dir{-};
(123,0); (123,-24) **\dir{-};
(81,0); (81,5) **\dir{-};(81,-24); (81,-29) **\dir{-};
(105,0); (105,5) **\dir{-};(105,-24); (105,-29) **\dir{-};
(45,0); (45,5) **\dir{-};(45,-24); (45,-29) **\dir{-};
(137,9)*{U^1};
(137,-33)*{U^2};
(-4,-8)*{p^{-1}(C_n)};
(-4,-16)*{=H\widetilde{C}_n};
(81,-12)*{\widetilde{C}_n};
(45,-35)*{\begin{xy}
   (-\el,\el); (\el,\el) **\dir{-}; (-\el,\el); (-\el,-\el) **\dir{-};
   (\el,\el); (\el,-\el) **\dir{-}; (-\el,-\el); (\el,-\el) **\dir{-};
   (20,0) *{q^{-1}(\overline{C}_n)=K\widetilde{C}_n};
   (-2.5, 2); (-2, 2.5) **\dir{-}; (-2.5, 0); (0, 2.5) **\dir{-};
   (-2.5, -2); (2, 2.5) **\dir{-}; (-1, -2.5); (2.5, 1) **\dir{-};
   (1, -2.5); (2.5, -1) **\dir{-};
   \end{xy}
};
(95,-35)*{\begin{xy}
   (-\el,\el); (\el,\el) **\dir{-}; (-\el,\el); (-\el,-\el) **\dir{-};
   (\el,\el); (\el,-\el) **\dir{-}; (-\el,-\el); (\el,-\el) **\dir{-};
   (6,0) *{A};
   \end{xy}
};
(115,-35)*{\begin{xy}   
   (-\el,\el); (\el,\el) **\dir{-}; (-\el,\el); (-\el,-\el) **\dir{-};
   (\el,\el); (\el,-\el) **\dir{-}; (-\el,-\el); (\el,-\el) **\dir{-};
   (6,0) *{B};
   @={(-1.65,2.2),(-0.55,2.2),(0.55,2.2),(1.65,2.2),(-1.1,1.1),(0,1.1),
   (1.1,1.1),(-1.65,0),(-0.55,0),(0.55,0),(1.65,0),(-1.1,-1.1),(0,-1.1),
   (1.1,-1.1),(-1.65,-2.2),(-0.55,-2.2),(0.55,-2.2),(1.65,-2.2)}
   @@{*{\cdot}}; 
   \end{xy}
};
(13.5, -2); (15.5, 0) **\dir{-};(13.5, -6); (19.5, 0) **\dir{-};
(13.5, -10); (23.5, 0) **\dir{-};(13.5, -14); (27, -0.5) **\dir{-};
(13.5, -18); (27, -4.5) **\dir{-};(13.5, -22); (27, -8.5) **\dir{-};
(15.5, -24); (27, -12.5) **\dir{-};(19.5, -24); (27, -16.5) **\dir{-};
(23.5, -24); (27, -20.5) **\dir{-};
(63, -2); (65, 0) **\dir{-};(63, -6); (69, 0) **\dir{-};
(63, -10); (73, 0) **\dir{-};(63, -14); (77, 0) **\dir{-};
(63, -18); (81, 0) **\dir{-};(63, -22); (85, 0) **\dir{-};
(77, -24); (87, -14) **\dir{-};(65, -24); (77.5, -11.5) **\dir{-};
(81,-8); (87,-2) **\dir{-};(69, -24); (77.5, -15.5) **\dir{-};
(83,-10); (87,-6) **\dir{-};(73, -24); (81, -16) **\dir{-};
(85,-12); (87,-10) **\dir{-};
(81, -24); (87, -18) **\dir{-};(85, -24); (87, -22) **\dir{-};
(123, -2); (125, 0) **\dir{-};(123, -6); (129, 0) **\dir{-};
(123, -10); (133, 0) **\dir{-};(123, -14); (136.5, -0.5) **\dir{-};
(123, -18); (136.5, -4.5) **\dir{-};(123, -22); (136.5, -8.5) **\dir{-};
(125, -24); (136.5, -12.5) **\dir{-};(129, -24); (136.5, -16.5) **\dir{-};
(133, -24); (136.5, -20.5) **\dir{-};
@={(17,-2),(21,-2),(25,-2),(29,-2),(33,-2),(37,-2),(19,-6),(23,-6),(27,-6),
(31,-6),(35,-6),(17,-10),(21,-10),(25,-10),(29,-10),(33,-10),(37,-10),(19,-14),
(23,-14),(27,-14),(31,-14),(35,-14),(17,-18),(21,-18),(25,-18),(29,-18),
(33,-18),(37,-18),(19,-22),(23,-22),(27,-22),(31,-22),(35,-22)}
@@{*{\cdot}}; 
@={(53,-2),(57,-2),(61,-2),(65,-2),(69,-2),(73,-2),(77,-2),(81,-2),(85,-2),
(89,-2),(93,-2),(97,-2),(55,-6),(59,-6),(63,-6),(67,-6),(71,-6),(75,-6),
(79,-6),(83,-6),(87,-6),(91,-6),(95,-6),(53,-10),(57,-10),(61,-10),(65,-10),
(69,-10),(77,-10),(85,-10),(89,-10),(93,-10),(97,-10),(55,-14),
(59,-14),(63,-14),(67,-14),(71,-14),(75,-14),(87,-14),(91,-14),
(95,-14),(53,-18),(57,-18),(61,-18),(65,-18),(69,-18),(73,-18),(77,-18),
(81,-18),(85,-18),(89,-18),(93,-18),(97,-18),(55,-22),(59,-22),(63,-22),
(67,-22),(71,-22),(75,-22),(79,-22),(83,-22),(87,-22),(91,-22),(95,-22)}
@@{*{\cdot}};
@={(113,-2),(117,-2),(121,-2),(125,-2),(129,-2),(133,-2),(115,-6),(119,-6),
(123,-6),(127,-6),(131,-6),(113,-10),(117,-10),(121,-10),(125,-10),(129,-10),
(133,-10),(115,-14),(119,-14),(123,-14),(127,-14),(131,-14),(113,-18),
(117,-18),(121,-18),(125,-18),(129,-18),(133,-18),(115,-22),(119,-22),
(123,-22),(127,-22),(131,-22)}
@@{*{\cdot}};

\end{xy}
} 

\caption{The decomposition of $N(\pcn)$ from Claim \ref{mainclaim}.}
\label{fig1}
\end{figure}

$\pcn = H\Cnt$ and $\qcn = K\Cnt$ (since $p(\Cnt)=\Cn$ and $q(\Cnt) = \Cnb$).
Let $A=\bigcup\left\{h\Cnt:h\in H \right.$ and $\left.h\Cnt\cap K\Cnt =
\emptyset\right\}$, $B=\bigcup\left\{h\Cnt:h\in H\right.$ and $\left.h\Cnt\cap
K\Cnt \neq \emptyset\right\}$.  The proof of the next claim contains the main
argument of the proof of the theorem (compare with Figure \ref{fig1}).

\begin{claim}
\label{mainclaim}
i) $A \neq \emptyset$.  ${\hat U_n} := A \;\cup\; U_n \;\cup\; \cup\{$all
bridge edges of $\pcn$ that intersect $A$ but not $\qcn\}$ is path
connected in $Y \cc \qcn$ and $\L$-unbounded. \\
ii) $B \;\cup\; \cup\{\mbox{all bridge edges of } \pcn \mbox{ that intersect }
B\}$ is $\L$-bounded.
\end{claim}

\begin{proof} If $h\Cnt \cap k\Cnt \neq \emptyset$ for some $h \in H, k \in K$,
  then $k^{-1}h\Cnt\cap\Cnt \neq \emptyset$, so letting $S = \{g\in H: g\Cnt
  \cap \Cnt \neq \emptyset\}$ we get $k^{-1}h =: s_0\in S$.  Now let $T$ be a
  set of right coset representatives for $K$ in $H$.  For each $s\in S$ write
  $s = k_s t_s$, where $t_s\in T$ and $k_s\in K$.  Then $k^{-1}h = k_{s_0}
  t_{s_0}$ and $h = (kk_{s_0}) t_{s_0}$, so if $h\Cnt \cap K\Cnt \neq
  \emptyset$ then $h$ lies in one of finitely many right cosets of $K$ in
  $H$. \\
  i) By assumption $[H:K] = \infty$, so there are in fact infinitely many $h
  \in H$ with $h\Cnt\cap K\Cnt = \emptyset$.  So $A$ is non-empty.  Recall that
  $\Cnt = \Cnh \cap \pcn$, so bridge edges of $\pcn$ contained
  in $\Cnh$ intersect $\pcn$ in $\Cnt$.  So if $h\Cnt \subseteq A$ then the
  bridge edges of $\pcn$ contained in $h\Cnh$ only intersect $\pcn$ in
  $h\Cnt\subseteq A$ and not in $\qcn$.  Then for $h\Cnt\subseteq A$ we get
  $h\Cnh \cup U_n = h\Cnt \cup \{$bridge edges of $\pcn$ in $h\Cnh \} \cup U_n
  \subseteq {\hat U_n}$, and $h\Cnh \cup U_n$ is path connected by $(\dagger)$,
  so ${\hat U_n}$ is path connected.  $(A\cup U_n) \cap \qcn = \emptyset$ by
  definition of $A$ and $U_n$, so ${\hat U_n}$ is path connected in $Y \cc
  \qcn$.  ${\hat U_n}$ is $\L$-unbounded since it contains $U_n$ which is
  $\K$-unbounded ($\K$-unbounded implies $\L$-unbounded).  \\
  ii) First note that $q(B)$ is finite since $q(K\Cnt)=\Cnb$ is finite and the
  other edges of $B$ are within a bounded edge distance of $K\Cnt$ (and $\Xb$
  is locally finite).  But then the union of {\em all} edges intersecting
  $B$ has finite image.
\end{proof}

Now write $Y \cc \qcn = {\bar U_n} \sqcup {\bar D_n}$, where ${\bar U_n}$ is
the path component containing ${\hat U_n}$ and ${\bar D_n}$ is the union of the
other components (possibly empty).  ${\bar U_n}$ is $\L$-unbounded since ${\hat
  U_n}$ is.  We can write $Y$ as $Y = \pcn \;\cup\; \{$bridge edges of $\pcn\}
\;\cup\; U_n$, and $A \;\cup\; B = \pcn$, so $Y = (A \;\cup\; U_n \;\cup\;
\{$bridge edges intersecting $A$ but not $\qcn\}) \;\cup\; (B \;\cup\;
\{$bridge edges intersecting $B\})$, and ${\bar U_n}$ contains the first
parenthetical term, so ${\bar D_n}$ must be contained in the second
parenthetical term, which by Claim $\ref{mainclaim}.ii)$ is $\L$-bounded.  So
for each $n$, $Y \cc \qcn = {\bar U_n} \sqcup {\bar D_n}$ where ${\bar U_n}$ is
$\L$-unbounded and ${\bar D_n}$ is $\L$-bounded.

Finally as in the proof of Lemma \ref{unlem} we can adjust the $\{\Cnb\}$
filtration so that $\qcn$ is a single $\L$-unbounded component for each
$n$, from which it follows that $\underset n \varprojlim\left\{\pi_0\left(Y\cc
\qcn\right)\right\}$ has cardinality $1$.  i.e. $\e{Y}{q} = 1$.
\endproof

\section{Construction of compatible filtrations}
\label{filts}

We now undertake the proof of Lemma \ref{compat}.  Recall
that we are working under the hypotheses of Theorem \ref{mt} with $Y$ a locally
finite graph, $p:Y\ra X$, $q:Y\ra\Xb$, $r:\Xb\ra X$, and $0<\e{Y}{p}=e<\infty$.
We restate Lemma \ref{compat}:\newline

\vspace{-6pt}

\noindent\textbf{Lemma 7.}{\em
  \ \ There exists a finite filtration $\{\Cnb\}$ of $\Xb$ and a finite
  filtration $\{\Cn\}$ of $X$ with $r(\Cnb)=\Cn$ and $Y \cc \pcn$ the disjoint
  union of $e$ $\{\pcn\}$-unbounded path components.  There are also finite
  subgraphs $\Cnt\subseteq \pcn$ and ${\hat \Cn} \subseteq N(\pcn)$ in $Y$ with
  ${\hat \Cn}\cap \pcn = \Cnt$, $q(\Cnt) = \Cnb$, $p(\Cnt) = \Cn$, and ${\hat
    \Cn}
  \cup U_n$ path connected.}\\
\endproof

The proof will start from the following general lemma.

\begin{lem}
\label{unlem}
Let $A$ be a locally finite graph and $f:A\ra B$ a covering map.  If $\e{A}{f}
= e$, $0<e<\infty$, then there exists a finite filtration $\{D_n\}$ of $B$ such
that $A \cc \pdn$ is the union of $e$ $\{\pdn\}$-unbounded path components for
all $n$.
\end{lem}
\begin{proof} Let $\{K_n\}$ be any finite filtration of $B$ and let
  $\L=\{p^{-1}(K_n)\}$.  By Corollary 4.7.5 of \cite{Ge} we may assume $A \cc
  p^{-1}(K_n)$ has exactly $e$ $\L$-unbounded components for all $n$.  Now for
  each $n$ let $D_n = K_n \bigcup \cup\{\mbox{finite components of }
  B\cc K_n\}$.  The union is finite (since $B$ is locally finite), so $D_n$ is 
  finite.  The image of a $\{\pdn\}$-bounded component of $A\cc \pdn$ would
  have to be contained in a finite component of $B\cc D_n$ (see the proof of
  Lemma $\ref{lft}.iv)$), so $A\cc \pdn$ consists exactly of $e$
  $\{\pdn\}$-unbounded components.
\end{proof}

\begin{rem} If $e=\infty$ in the hypotheses of Lemma \ref{unlem}, then there
  exists a finite filtration $\{\Cn\}$ of $X$ such that $Y\cc \pcn$ equals the
  union of $\{\pcn\}$-unbounded path components, with the number of
  path components an increasing function of $n$.
\end{rem}

{\bf PROOF of Lemma \ref{compat}.}  Let $\{\Cnb\}$ be a finite filtration of
$\Xb$.  Define $\Cn = r(\Cnb)$.  $\{\Cn\}$ is a finite filtration of $X$.  By
the (proof of) Lemma \ref{unlem} (after renumbering) for each $n$ there exists
a finite subgraph $F_n$ of $X$ with the property that $\{C_n \cup F_n\}$ is a
finite filtration of $X$ and $Y \cc p^{-1}(C_n \cup F_n)$ is the union of $e$
$\{p^{-1}(C_n \cup F_n)\}$-unbounded path components.  Now for each $n$ choose
a finite subgraph ${\bar F_n}$ in $\Xb$ with $r({\bar F_n})=F_n$, choose a
subsequence filtration of $\{\Cnb \cup {\bar F_n}\}$, and then rename the $\Xb$
subsequence filtration to be $\{\Cnb\}$ and the $X$ filtration to be $\{\Cn =
r(\Cnb)\}$, a subsequence of the old $\{C_n \cup F_n\}$ filtration.  To
summarize: we have a finite filtration $\{\Cnb\}$ of $\Xb$ and a finite
filtration $\{\Cn\}$ of $X$ with the property that $r(\Cnb) = \Cn$ and $Y \cc
\pcn$ is the union of $e$ $\{\pcn\}$-unbounded path components for all $n$.

Let $\K = \{\pcn\}$.

Choose a finite subgraph $\Cnt$ in $Y$ with $q(\Cnt) = \Cnb$.  The given
lifts $\Cnt$ will not in general satisfy the properties we will require, but we
can use them to construct all new filtrations that will, which we do now.  Fix
$n$.  Let $U^1, \ldots, U^e$ be the $\K$-unbounded components of $Y \cc \pcn$
and define $U_n = U^1 \cup \cdots \cup U^e$.  Our next aim is to define new
$\Cnt$ so that $U_n\cup \Cnt$ is path connected (up to bridge edges of $\pcn$).
Claim \ref{hat} tells us how to adjust our original choices (compare with
Figure \ref{fig2}).
  
\begin{claim}
\label{hat}
There exists a finite subgraph ${\hat C_n}$ of $N(\pcn)$ containing $\Cnt$
such that ${\hat C_n} \cup U_n$ is path connected.
\end{claim}
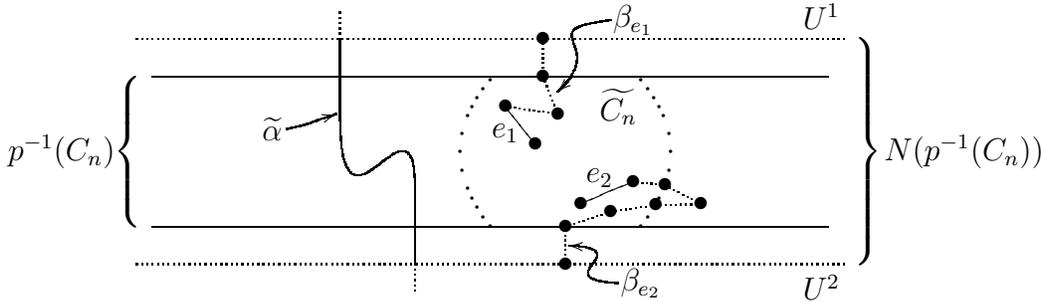
\begin{figure}[h]
{\fontsize{12}{12}\selectfont
\begin{xy}
(28,5); (122,5) **\dir{.};
(28,-25); (122,-25) **\dir{.};
(30,0); (120,0) **\dir{-};
(30,-20); (120,-20) **\dir{-};
(119,8) *{U^1};
(119,-28) *{U^2};
(27,-10)*{\left\{
   \begin{xy}
   (0,10.5)*{};
   (0,-10.5)*{};
   \end{xy}
   \right.
};
(18,-10) *{p^{-1}(C_n)};
(125,-10) *{\left.
   \begin{xy}
   (0,15)*{};
   (0,-15)*{};
   \end{xy}
   \right\}
};
(138,-10) *{N(p^{-1}(C_n))};
(75,-.5); (75,-20) **\crv{~*{\cdot}(67,-10)};
(95,-.5); (95,-20) **\crv{~*{\cdot}(103,-10)}
?!{(94,-14); (102,-15)}="i1"
?!{(102,-16); (93,-18)}="i2";
(92,-4) *{\widetilde{C_n}};
@={(81,-9),(77,-4),(84,-5),(82,0),(82,5)}
@@{*{\bullet}};
(81,-9); (77,-4) **\dir{-} ?+(-2,-1)*{e_1};
(77,-4); (84,-5) **\dir{.};
(82,0) **\dir{.} ?="c";
(82,0); (82,5) **\dir{.};
"c"+(1,.5) ; (91,7.5) **\crv{(91,4.5) & (80,6)} ?<*\dir{<}
?>+(3,0)*{\beta_{e_1}};
@={(87,-17),(94,-14),"i1",(103,-17),"i2",(91,-18),(85,-20),(85,-25)}
@@{*{\bullet}};
(87,-17); (94,-14) **\dir{-} ?+(-1,2)*{e_2};
(94,-14); "i1" **\dir{.};
(103,-17) **\dir{.};
"i2" **\dir{.};
(91,-18) **\dir{.};
(85,-20) **\dir{.};
(85,-25) **\dir{.};
(85.5,-22.5); (92,-27.5)  **\crv{(95,-23) & (82,-27)} ?<*\dir{<}
?>+(3,0)*{\beta_{e_2}};
(55,5); (55,8.5) **\dir{.};
(55,5); (55,-5) **\dir{-};
(65,-25); (65,-28.5) **\dir{.};
(65,-25); (65,-15) **\dir{-};
(55,-5); (65,-15) **\crv{(55,-25) & (65,0)};
(48,-7); (54.5,-5) **\crv{(52,-6) & (50,-6.5)} ?>*\dir{>}
?<+(-2,0)*{\widetilde{\alpha}};
\end{xy}
} 
\caption{The construction of ${\hat C_n}$ in Claim \ref{hat}.}
\label{fig2}
\end{figure}
\begin{proof} Let $\alpha$ be a finite edge path (i.e. union of edges) in $Y$
  starting in $U_n$ with the property that $\alpha$ intersects each $U^i$
  non-trivially.  Let ${\tilde \alpha}$ be the edges of $\alpha$ contained in
  $N(\pcn)$.  Then ${\tilde \alpha} \cup U_n$ is path connected.  Now for each
  cell $e$ of $\Cnt$ choose a finite edge path $\beta_e$ in $N(\pcn)$ to $U_n$.
  Take $\displaystyle {\hat C_n} = \Cnt \cup {\tilde \alpha} \cup \left(
    \bigcup_{e\subset \Cnt} \beta_e \right)$.
\end{proof}

Let $\Cnt' = {\hat \Cn} \cap \pcn \supseteq \Cnt$, and let $\Cnb' = q(\Cnt')$
($p(\Cnt')=\Cn$, so we do not rename for $X$).  $\{\Cnb'\}$ may no longer be a
chain of subsets, but since we only added finitely many cells to each term of
the original filtration we can find a subsequence which is a finite filtration.
Now renumber and then rename all filtrations to get finite filtrations
$\{\Cnb\}$ for $\Xb$, $\{\Cn\}$ for $X$, and a finite lift $\Cnt$ for each
$\Cnb$ with the property that there is a finite ${\hat \Cn} \subseteq N(\pcn)$
containing $\Cnt$ with $U_n \cup {\hat \Cn}$ path connected -- note that we
only took a subsequence of our previous $X$ filtration $\{\Cn\}$, so we still
have $Y \cc \pcn = U_n$, where $U_n = U^1 \cup \cdots \cup U^e$ is the same
union of $\K$-unbounded path components as before (of course numbering and the
meaning of $\K$ have changed, but the components and their $\K$-unboundedness
have not).
\endproof 

\section{Filtered ends of pairs of groups}
\label{corollary}

We now apply this to groups.  The general situation will be as follows: \\
{\fontsize{11}{8}\selectfont
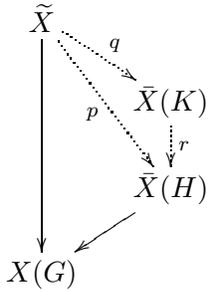
\xymatrix@=15pt{ {\widetilde X} \ar@{.>}[dr]^q \ar@{.>}[ddr]_p \ar[ddd] & \\
 & {\bar X}(K) \ar@{.>}[d]^r \\
 & {\bar X}(H) \ar[dl] \\
 X(G) &  }
} 
\\ \\
\noindent
$G$ is finitely generated infinite, $K\leq H \leq G$, and $X(G)$
has finite $1$-skeleton and fundamental group $G$ (we write ``$\Xb(L)$'' for
the covering space of $X$ corresponding to $L\leq G$).  We apply Theorem
\ref{mt} to the dotted triple of covering spaces to obtain Corollary \ref{cor1}
on filtered ends of pairs of groups.

{\bf PROOF of Corollary \ref{cor1}.}
  The case $[H:K]<\infty$ is standard (\cite{KR}, $2.4.v)$); the
  case $[H:K]=\infty$ is the statement of Theorem \ref{mt}.
\endproof

Corollary \ref{cor1} and its contrapositive give the following picture of
filtered end behavior above and below a subgroup $H$ with $\e{G}{H} = n$:\\

\def\l1{24} 
\def\m2{12} 
\def\n3{8}  
\def\o32{16}
\def\w1{13} 
\def\s1{20} 
\def\t1{40} 
\def\u1{72} 
\def\v1{57} 

{\fontsize{10}{12}\selectfont
\begin{xy}
(0,\l1)*+{G}="G";
(0,0)*+{H}="H";
(-\w1,\o32)*{}="L1";
(-\w1,\n3)*{}="L2";
(\w1,\o32)*{}="R1";
(\w1,\n3)*{}="R2";
"G";"L1" **\crv{(-\w1,\l1)};
"G";"R1" **\crv{(\w1,\l1)};
"H";"L2" **\crv{(-\w1,0)};
"H";"R2" **\crv{(\w1,0)};
"L1";"L2" **\dir{-};
"R1";"R2" **\dir{-};
"L1";"R1" **\dir{-};
(15,\o32)*{}; (\t1,\o32)*{} **\dir{.};
"L2";"R2" **\dir{-};
(15,0)*{}; (\t1,0)*{} **\dir{.};
(0,20)*{e(G,K)=0};
(\s1,20)*[r]{[G:K]<\infty};
(0,\m2)*{e(G,K)=\infty};
(\s1,\m2)*[r]{[K:H]=\infty};
(0,4)*{e(G,K)=n};
(\s1,4)*[r]{[K:H]<\infty};
(0,-\o32)*+{1}="1";
(-\w1,-\n3)*{}="M";
(\w1,-\n3)*{}="N";
"H";"M" **\crv{(-\w1,0)};
"H";"N" **\crv{(\w1,0)};
"1";"M" **\crv{(-\w1,-\o32)};
"1";"N" **\crv{(\w1,-\o32)};
"M";"N" **\dir{-};
(15,-\n3)*{}; (\t1,-\n3)*{} **\dir{.};
(0,-5)*{e(G,K)=n};
(\s1,-5)*[r]{[H:K]<\infty};
(0,-12)*{e(G,K)=1};
(\s1,-12)*[r]{[H:K]=\infty};
(15,8)*{}; (\u1,8)*{} **\dir{.};
(9,\l1)*{}; (\u1,\l1)*{} **\dir{.};
(9,-\o32)*{}; (\u1,-\o32)*{} **\dir{.};
(\v1,\o32)*+{1<\e{G}{H}<\infty}="S1";
(\v1,-5)*+{0<\e{G}{H}<\infty}="T1";
{\ar "S1"; (\v1,\l1)};
{\ar "S1"; (\v1,8)};
{\ar "T1"; (\v1,8)};
{\ar "T1"; (\v1,-\o32)};
\end{xy}
} 

\end{document}